\magnification 1200
\input plainenc
\input amssym
\fontencoding{T2A}
\inputencoding{cp1251}
\tolerance 4000
\relpenalty 10000
\binoppenalty 10000
\parindent 1.5em

\hsize 17truecm
\vsize 24truecm
\hoffset 0truecm
\voffset -1truecm

\font\TITLE labx1440
\font\tenrm larm1000
\font\cmtenrm cmr10
\font\tenit lati1000
\font\tenbf labx1000
\font\teni cmmi10 \skewchar\teni '177
\font\tensy cmsy10 \skewchar\tensy '60
\font\tenex cmex10
\font\teneufm eufm10
\font\eightrm larm0800
\font\cmeightrm cmr8
\font\eightit lati0800
\font\eightbf labx0800
\font\eighti cmmi8 \skewchar\eighti '177
\font\eightsy cmsy8 \skewchar\eightsy '60
\font\eightex cmex8
\font\eighteufm eufm8

\font\cmsixrm cmr6

\font\sixbf labx0600
\font\sixi cmmi6 \skewchar\sixi '177
\font\sixsy cmsy6 \skewchar\sixsy '60
\font\sixeufm eufm6

\font\cmfiverm cmr5

\font\fivebf labx0500
\font\fivei cmmi5 \skewchar\fivei '177
\font\fivesy cmsy5 \skewchar\fivesy '60
\font\fiveeufm eufm5
\font\tencmmib cmmib10 \skewchar\tencmmib '177
\font\eightcmmib cmmib8 \skewchar\eightcmmib '177
\font\sevencmmib cmmib7 \skewchar\sevencmmib '177
\font\sixcmmib cmmib6 \skewchar\sixcmmib '177
\font\fivecmmib cmmib5 \skewchar\fivecmmib '177
\newfam\cmmibfam
\textfont\cmmibfam\tencmmib \scriptfont\cmmibfam\sevencmmib
\scriptscriptfont\cmmibfam\fivecmmib
\def\tenpoint{\def\rm{\fam0\tenrm}\def\it{\fam\itfam\tenit}%
	\def\bf{\fam\bffam\tenbf}
	\textfont0\cmtenrm \scriptfont0\cmsevenrm \scriptscriptfont0\cmfiverm
  	\textfont1\teni \scriptfont1\seveni \scriptscriptfont1\fivei
  	\textfont2\tensy \scriptfont2\sevensy \scriptscriptfont2\fivesy
  	\textfont3\tenex \scriptfont3\tenex \scriptscriptfont3\tenex
  	\textfont\itfam\tenit
	\textfont\bffam\tenbf \scriptfont\bffam\sevenbf
	\scriptscriptfont\bffam\fivebf
	\textfont\eufmfam\teneufm \scriptfont\eufmfam\seveneufm
	\scriptscriptfont\eufmfam\fiveeufm
	\textfont\cmmibfam\tencmmib \scriptfont\cmmibfam\sevencmmib
	\scriptscriptfont\cmmibfam\fivecmmib
	\normalbaselineskip 12pt
	\setbox\strutbox\hbox{\vrule height8.5pt depth3.5pt width0pt}%
	\normalbaselines\rm}
\def\eightpoint{\def\rm{\fam 0\eightrm}\def\it{\fam\itfam\eightit}%
	\def\bf{\fam\bffam\eightbf}%
	\textfont0\cmeightrm \scriptfont0\cmsixrm \scriptscriptfont0\cmfiverm
	\textfont1\eighti \scriptfont1\sixi \scriptscriptfont1\fivei
	\textfont2\eightsy \scriptfont2\sixsy \scriptscriptfont2\fivesy
	\textfont3\eightex \scriptfont3\eightex \scriptscriptfont3\eightex
	\textfont\itfam\eightit
	\textfont\bffam\eightbf \scriptfont\bffam\sixbf
	\scriptscriptfont\bffam\fivebf
	\textfont\eufmfam\eighteufm \scriptfont\eufmfam\sixeufm
	\scriptscriptfont\eufmfam\fiveeufm
	\textfont\cmmibfam\eightcmmib \scriptfont\cmmibfam\sixcmmib
	\scriptscriptfont\cmmibfam\fivecmmib
	\normalbaselineskip 11pt
	\abovedisplayskip 5pt
	\belowdisplayskip 5pt
	\setbox\strutbox\hbox{\vrule height7pt depth2pt width0pt}%
	\normalbaselines\rm
}

\def\No{\char 157}
\def\empty{}

\catcode`\@ 11
\catcode`\" 13
\def"#1{\ifx#1<\char 190\relax\else\ifx#1>\char 191\relax\else #1\fi\fi}

\def\newl@bel#1#2{\expandafter\def\csname l@#1\endcsname{#2}}
\openin 11\jobname .aux
\ifeof 11
	\closein 11\relax
\else
	\closein 11
	\input \jobname .aux
	\relax
\fi

\newcount\c@section
\newcount\c@subsection
\newcount\c@subsubsection
\newcount\c@equation
\newcount\c@bibl
\c@section 0
\c@subsection 0
\c@subsubsection 0
\c@equation 0
\c@bibl 0
\def\lab@l{}
\def\label#1{\immediate\write 11{\string\newl@bel{#1}{\lab@l}}%
	\ifhmode\unskip\fi}
\def\eqlabel#1{\rlap{$(\equation)$}\label{#1}}

\def\section#1{\global\advance\c@section 1
	{\par\vskip 3ex plus 0.5ex minus 0.1ex
	\rightskip 0pt plus 1fill\leftskip 0pt plus 1fill\noindent
	{\bf\S\thinspace\number\c@section .~#1}\par\penalty 25000%
	\vskip 1ex plus 0.25ex}
	\gdef\lab@l{\number\c@section.}
	\c@subsection 0
	\c@subsubsection 0
	\c@equation 0
}
\def\subsection{\global\advance\c@subsection 1
	\par\vskip 1ex plus 0.1ex minus 0.05ex{\bf\number\c@subsection. }%
	\gdef\lab@l{\number\c@section.\number\c@subsection}%
	\c@subsubsection 0\c@equation 0%
}
\def\subsubsection{\global\advance\c@subsubsection 1
	\par\vskip 1ex plus 0.1ex minus 0.05ex%
	{\bf\number\c@subsection.\number\c@subsubsection. }%
	\gdef\lab@l{\number\c@section.\number\c@subsection.%
		\number\c@subsubsection}%
}
\def\equation{\global\advance\c@equation 1
	\gdef\lab@l{\number\c@section.\number\c@subsection.%
	\number\c@equation}{\rm\number\c@equation}
}
\def\bibitem#1{\global\advance\c@bibl 1
	\par[\number\c@bibl]%
	\gdef\lab@l{\number\c@bibl}\label{#1}
}
\def\ref@ref#1.#2:{\def\REF@{#2}\ifx\REF@\empty{\S\thinspace#1}%
	\else\ifnum #1=\c@section {#2}\else {\S\thinspace#1.#2}\fi\fi
}
\def\ref@eqref#1.#2.#3:{\ifnum #1=\c@section\ifnum #2=\c@subsection
	{(#3)}\else{#2\thinspace(#3)}\fi\else{\S\thinspace#1.#2\thinspace(#3)}\fi
}
\def\ref#1{\expandafter\ifx\csname l@#1\endcsname\relax
	{\bf ??}\else\edef\mur@{\csname l@#1\endcsname :}%
	{\expandafter\ref@ref\mur@}\fi
}
\def\eqref#1{\expandafter\ifx\csname l@#1\endcsname\relax
	{(\bf ??)}\else\edef\mur@{\csname l@#1\endcsname :}%
	{\expandafter\ref@eqref\mur@}\fi
}
\def\cite#1{\expandafter\ifx\csname l@#1\endcsname\relax
	{\bf ??}\else\hbox{\bf\csname l@#1\endcsname}\fi
}

\def\superalign#1{\openup 1\jot\tabskip 0pt\halign to\displaywidth{%
	\tabskip 0pt plus 1fil$\displaystyle ##$&%
	\tabskip 0pt\hss $\displaystyle ##{}$&%
	$\displaystyle {}##$\hss\tabskip 0pt plus 1fil&%
	\tabskip 0pt\hss ##\crcr #1\crcr}}

\def\fB{\frak B}

\def\fD{\frak D}
\def\fH{\frak H}
\def\Re{\mathop{\rm Re}}

\def\dom{\mathop{\rm dom}}
\def\ind{\mathop{\rm ind}}
\catcode`\" 12
\def\bolddelta{\mathchar"0\hexnumber@\cmmibfam 0E}
\catcode`\" 13
\def\Wo{{\mathpalette\Wo@{}}W}
\def\Wo@#1{\setbox0\hbox{$#1 W$}\dimen@\ht0\dimen@ii\wd0\raise0.65\dimen@%
\rlap{\kern0.35\dimen@ii$#1{}^\circ$}}
\catcode`\@ 12

\def\proof{\par\medskip{\rm Д$\,$о$\,$к$\,$а$\,$з$\,$а$\,$т$\,$е$\,$л$\,$ь%
	$\,$с$\,$т$\,$в$\,$о.}\ }
\def\endproof{{\parfillskip 0pt\hfill$\square$\par}\medskip}

\immediate\openout 11\jobname.aux


\frenchspacing\rm
\ \par\vskip 1truecm
\leftline{\rm УДК~517.984}
{\par\medskip\noindent\rightskip=0pt plus 1fill\leftskip=0pt plus 1fill
\noindent
\TITLE Теоремы о представлении и вариационные принципы для самосопряжённых
операторных матриц\par\bigskip\rm А.$\,$А.~Владимиров\par}

\rm
\vskip 1truecm
\section{Введение}\label{par:1}
\subsection\label{pt:1.1}
Рассмотрим следующую ситуацию. Пусть дано некоторое гильбертово простран\-ство $\fH$,
разложенное в ортогональную прямую сумму $\fH_1\oplus\fH_2$ двух своих замкнутых
подпространств. Пусть при этом дополнительно фиксированы четыре оператора
$T^\circ_{11}\colon\dom T^\circ_{11}\to\fH_1$, $T_{12}^\circ\colon\dom T^\circ_{22}
\to\fH_1$, $T^\circ_{21}\colon\dom T^\circ_{11}\to\fH_2$ и $T_{22}^{\circ}\colon
\dom T_{22}^\circ\to\fH_2$, задающие симметрическую операторную матрицу
$$
	T^{\circ}\rightleftharpoons\pmatrix{T_{11}^{\circ}& T_{12}^{\circ}\cr
		T_{21}^{\circ}&T_{22}^{\circ}\cr}\leqno(\equation)
$$\label{eq:1.99}%
с областью определения $\dom T_{11}^\circ\oplus\dom T_{22}^\circ$. Наконец, пусть
симметрический оператор $T_{11}^{\circ}$ ограничен снизу, а симметрический
оператор $T_{22}^{\circ}$~--- сверху.

В работе [\cite{KLT}] для трёх различных типов взаимных соотношений между
операторами $T_{11}^{\circ}$, $T_{22}^{\circ}$ и $T_{21}^{\circ}$ была развита
вариационная техника, позволяющая исследовать свойства некоторых частей дискретного
спектра замыкания оператора $T^\circ$. Основной целью настоящей статьи является
указание ряда общих фактов, частными случаями которых выступают результаты работы
[\cite{KLT}], а также некоторые аналогичные им.

\subsection
Основу предлагаемого нами подхода составляет теория {\it оснащённых пространств\/}
(см., например, [\cite{LM}, Глава~1, п.~2.4] или [\cite{BSh:1983},
Добавление~1, \S\S$\;$2.3]). В настоящей статье под оснащённым пространством будет
пониматься совокупность ${\bf A}\rightleftharpoons\{\fB,\fD^+,\fD^-,I^+,I^-\}$
из трёх банаховых пространств $\fB$, $\fD^+$ и $\fD^-$ и двух инъективных вложений
$I^+\colon\fD^+\to\fB$ и $I^-\colon\fB\to\fD^-$ с плотными образами:
$$
	\fD^+\buildrel I^+\over\to\fB\buildrel I^-\over\to\fD^-.
$$
{\it Операторами в оснащении $\bf A$\/} мы далее будем сокращённо называть
ограниченные операторы, отображающие пространство $\fD^+$ в пространство $\fD^-$.
Класс операторов в оснащении $\bf A$ мы будем обозначать символом
${\cal B}({\bf A})$.

\subsection
Структура статьи имеет следующий вид. В~\ref{par:20} устанавливаются используемые
далее теоремы о представлении замкнутого неограниченного оператора в банаховом
пространстве линейным пучком ограниченных операторов в оснащении. Эти теоремы
содержат как простые частные случаи классические результаты о представлении
секториального оператора полуторалинейной формой [\cite{Kato:1972},
Гл.~\hbox{VI}, \S$\;$2]. В~\ref{par:2} вводится процедура {\it углового
расширения\/} симметрической операторной матрицы, опирающаяся
на результаты~\ref{par:20} и тесно связанная со стандартной процедурой расширения
(и псевдорасширения) оператора по Фридрихсу. В~\ref{par:3} на основе вариационных
принципов для самосопряжённых оператор-функций устанавливаются вариационные
принципы для угловых расширений. Наконец, в~\ref{par:4} обсуждаются некоторые
дополнительные применения результатов из~\ref{par:20}.


\section{Теоремы о представлении}\label{par:20}
\subsection\label{par:20.1}
Пусть зафиксированы некоторое оснащение ${\bf A}\rightleftharpoons\{\fB,\fD^+,
\fD^-,I^+,I^-\}$ и связанный с ним оператор $T\in{\cal B}({\bf A})$. Символом
$T^\bullet$ мы далее будем обозначать оператор $(I^-)^{-1}T(I^+)^{-1}$,
действующий в пространстве $\frak B$ и, вообще говоря, неограниченный.

\subsubsection\label{20:pred}
{\it Пусть при некотором $A\in{\cal B}(\fB)$ оператор $T-I^-AI^+$ обладает
ограниченным обратным. Тогда оператор $T^\bullet$ замкнут и плотно определён.
}

\proof
Легко видеть справедливость равенства
$$
	T^\bullet-A=(I^-)^{-1}\cdot (T-I^-AI^+)\,(I^+)^{-1}.\leqno(\equation)
$$\label{eq:2.60}%
Соответственно, оператор $T^\bullet-A$ является обратным к всюду определённому
ограниченному оператору $I^+\cdot(T-I^-AI^+)^{-1}I^-$, а потому замкнут.
Этот факт автоматически означает замкнутость оператора $T^\bullet$.

Далее, область определения оператора $T^\bullet$ совпадает с областью определения
оператора $T^\bullet-A$, а потому и с областью значений оператора $I^+\cdot
(T-I^-AI^+)^{-1}I^-$. Ввиду предполагаемой плотности образов операторов $I^\pm$
это означает плотность подмножества $\dom T^\bullet$ в пространстве $\fB$.
\endproof

\subsubsection\label{20:sod}
{\it Всякое значение $\lambda\in\Bbb C$, для которого оператор $T-\lambda I^-I^+$
обладает ограниченным обратным, принадлежит резольвентному множеству оператора
$T^\bullet$ и удовлетворяет равенству
$$
	(T^\bullet-\lambda)^{-1}=I^+\cdot(T-\lambda I^-I^+)^{-1}I^-.
$$
}

Справедливость данного предложения, по существу, установлена в ходе доказатель\-ства
предложения \ref{20:pred}.

\subsubsection\label{prop:2.2}
{\it Пусть при некотором $A\in{\cal B}(\frak B)$ оператор $T-I^-AI^+$ обладает
ограниченным обратным. Тогда для любого $\lambda\in\varrho(T^\bullet)$ оператор
$T-\lambda I^-I^+$ также обладает ограниченным обратным.
}

\proof
Введём сокращения $T^\natural(\lambda)\rightleftharpoons T-\lambda I^-I^+$,
$T^\natural(A)\rightleftharpoons T-I^- AI^+$, $R_\lambda\rightleftharpoons
(T^\bullet-\lambda)^{-1}$ и $R_A\rightleftharpoons (T^\bullet-A)^{-1}$, а также
обозначим через $S\colon\fB\to\fD^+$ ограниченный оператор
$$
	S\rightleftharpoons [T^\natural(A)]^{-1}I^-\cdot[1+(\lambda-A)R_\lambda].
	\leqno(\equation)
$$\label{eq:2.59}%
Заметим, что имеют место равенства
$$
	\belowdisplayskip=0pt
	\leqalignno{R_\lambda-R_A&=R_A\cdot[(T^\bullet-A)-(T^\bullet-\lambda)]
		R_\lambda\cr &=R_A\cdot (\lambda-A)R_\lambda,&\eqlabel{eq:2.58}}
$$
$$
	\abovedisplayskip=0pt
	\superalign{&T^\natural(\lambda)S&=[I^--I^-(\lambda-A)R_A]\cdot
		[1+(\lambda-A)R_\lambda]&[\eqref{eq:2.59}, \eqref{eq:2.60}]\cr
		&&=I^-\cdot[1+(\lambda-A)(R_\lambda-R_A)-(\lambda-A)R_A
		\cdot(\lambda-A)R_\lambda]\cr
		\eqlabel{eq:2.61}&&=I^-.&[\eqref{eq:2.58}]\cr
	}
$$
Соответственно, равенства
$$
	\superalign{&T^\natural(\lambda)\cdot [1+S(\lambda-A)I^+]
		[T^\natural(A)]^{-1}&=T^\natural(A)\cdot [T^\natural(A)]^{-1}
		&[\eqref{eq:2.61}]\cr &&=1\cr
	}
$$
означают существование у оператора $T^\natural(\lambda)$ ограниченного правого
обратного. Вытека\-ющая из ограниченной обратимости оператора $T^\bullet-\lambda$
инъективность оператора $T^\natural(\lambda)$ означает потому ограниченную
обратимость последнего.
\endproof

Сказанное означает, что если хотя бы для одного оператора $A\in{\cal B}(\fB)$
оператор $T-I^-AI^+$ обладает ограниченным обратным, то спектр оператора $T^\bullet$
в точности совпадает со спектром линейного пучка $T^\natural\colon\Bbb C\to
{\cal B}({\bf A})$ вида $T^\natural(\lambda)\rightleftharpoons T-\lambda I^-I^+$.
При этом также справедливо следующее предложение.

\subsubsection\label{20:ker}
{\it Независимо от выбора значения $\lambda\in\Bbb C$ выполняется равенство
$$
	\belowdisplayskip=0pt
	\ker(T^\bullet-\lambda)=I^+\ker T^\natural(\lambda).
$$
}

\proof
Любой вектор $x\in\ker(T^\bullet-\lambda)$ допускает представление в виде $x=I^+y$,
где вектор $y\in\fD^+$ удовлетворяет равенству $(I^-)^{-1}T^\natural(\lambda)y=0$.
Соответственно, имеет место вложение
$$
	\ker(T^\bullet-\lambda)\subseteq I^+\ker T^\natural(\lambda).
$$
С другой стороны, для любого вектора $y\in\ker T^\natural(\lambda)$ выполняется
равенство $Ty=I^-(\lambda I^+y)$, означающее принадлежность вектора
$x\rightleftharpoons I^+y$ области определения оператора $T^\bullet$, а тогда
и ядру оператора $T^\bullet-\lambda$.
\endproof

\subsection
Рассмотрим в качестве примера оснащение
$$
	{\bf A}_0\rightleftharpoons\{L_2[0,1],\tilde W_2^1[0,1],L_2[0,1],
		I,{\rm id}\}.
$$
Здесь использовано обозначение
$$
	\tilde W_2^1[0,1]\rightleftharpoons\{y\in W_2^1[0,1]\;:\;y(0)=y(1)\},
	\leqno(\equation)
$$\label{eq:Wtild}%
через $I\colon\tilde W_2^1[0,1]\to L_2[0,1]$ обозначен оператор вложения, а через
$\rm id$~--- тождественный оператор в пространстве $L_2[0,1]$. Свяжем с этим
оснащением оператор $T\in{\cal B}({\bf A}_0)$ вида $Ty\rightleftharpoons -iy'$.
Тогда оператор $T-I$ обладает обратным интегральным оператором с ядром
$$
	K(x,t)=\cases{\displaystyle{ie^{i(x-t)}\over 1-e^i}&при $x\geqslant t$,\cr
		\displaystyle{-ie^{i(x-t)}\over 1-e^{-i}}&при $x\leqslant t$,}
$$
что, согласно предложениям~\ref{20:pred} и \ref{prop:2.2}, означает замкнутость
и плотную определённость оператора $T^\bullet$, а также совпадение его спектра
со спектром пучка $T^\natural$.

Следующие примеры показывают, что условия из предложений~\ref{20:pred}
и~\ref{prop:2.2} носят содержательный характер.

Рассмотрим оснащение
$$
	{\bf A}_1\rightleftharpoons\{L_2[0,1],\tilde W_2^1[0,1],
		\tilde W_2^{-1}[0,1],I,I^*\},
$$
где через $\tilde W_2^{-1}[0,1]$ обозначено пространство, сопряжённое пространству
$\tilde W_2^1[0,1]$. В этом случае оператор $T\in{\cal B}({\bf A}_1)$ вида
$Ty\rightleftharpoons -iy'$ вполне непрерывен, что, ввиду полной непрерывности
вложения $I$, означает отсутствие для оператора $T-I^*AI$ ограниченного обратного
независимо от выбора оператора $A\in{\cal B}(L_2[0,1])$. Оператор $T^\bullet$ имеет
тот же вид, что и в рассмотренном выше случае оснащения ${\bf A}_0$, однако
спектром пучка $T^\natural$ в рассматриваемой ситуации выступает вся плоскость
$\Bbb C\neq\sigma(T^\bullet)$.

Рассмотрим оснащение
$$
	{\bf A}_2\rightleftharpoons\{L_2[0,1],\tilde W_2^2[0,1],L_2[0,1],
		J,{\rm id}\},
$$
где использовано обозначение
$$
	\tilde W_2^2[0,1]\rightleftharpoons\{y\in W_2^2[0,1]\;:\;y(0)=y(1),
		\; y'(0)=y'(1)\},
$$
а через $J\colon\tilde W_2^2[0,1]\to L_2[0,1]$ обозначен оператор вложения.
В этом случае оператор $T^\bullet$ не является замкнутым.

Наконец, рассмотрим оператор $T\in{\cal B}({\bf A}_1)$ вида
$$
	Ty\rightleftharpoons\sum_{n=1}^\infty 2^{-n}y(\zeta_n)
		\bolddelta_{\zeta_n},
$$
где $\{\zeta_n\}_{n=1}^\infty$~--- произвольно фиксированная плотная на отрезке
$[0,1]$ последовательность, а через $\bolddelta_\zeta\in\tilde W_2^{-1}[0,1]$
обозначена дельта-функция с сосредоточенным в точке $\zeta\in [0,1]$ носителем.
В этом случае область определения оператора $T^\bullet$ содержит лишь нулевой
вектор пространства $L_2[0,1]$.


\section{Угловые расширения}\label{par:2}
\subsection\label{par:2.1}
Классическая процедура расширения полуограниченного самосопряжённого оператора
по Фридрихсу на языке теории оснащённых пространств выглядит следующим образом.
Пусть $T^\circ$~--- действующий в гильбертовом пространстве $\fH$ симметрический
оператор (вообще говоря, неограниченный), и пусть число $\gamma\in\Bbb R$
удовлетворяет условию
$$
	(\forall y\in\dom T^\circ)\qquad\langle (T^\circ-\gamma)y,y\rangle_\fH
		\geqslant\|y\|_\fH^2.
$$
Обозначим через $\fD$ пополнение линейного множества $\dom T^\circ$ по норме
$$
	\|y\|_\fD\rightleftharpoons\langle (T^\circ-\gamma)y,y\rangle_\fH^{1/2},
$$
а через $I\colon\fD\to\fH$~--- соответствующий оператор вложения. Наконец, обозначим
через $\fD^*$ сопряжённое к $\fD$ пространство, после чего рассмотрим оснащение
${\bf A}\rightleftharpoons\{\fH,\fD,\fD^*,I,I^*\}$. Тривиальное тождество
$$
	(\forall y\in I^{-1}\dom T^\circ)\qquad\langle [I^* T^\circ I-
		\gamma I^*I]y,y\rangle_\fH=\|y\|^2_\fD\leqno(\equation)
$$\label{eq:3.99}%
означает, что оператор $I^*T^\circ I$ допускает однозначное продолжение
по непрерывности до некоторого оператора $T\in{\cal B}({\bf A})$. Далее следует
воспользоваться следующей теоремой.

\subsubsection\label{prop:St} [\cite{Kato:1972}, Гл.~V, Следствие~3.3]
{\it Пусть $\fD$~--- гильбертово пространство, и пусть числовая область значений
$$
	W(A)\rightleftharpoons\bigl\{\lambda\in\Bbb C\;:\;(\exists y\in\fD:
		\|y\|_{\fD}=1)\quad\langle Ay,y\rangle=\lambda\bigr\}
$$
ограниченного оператора $A\colon\fD\to\fD^*$ отделена от нуля. Тогда оператор $A$
обладает ограниченным обратным.
}

\medskip\noindent
С учётом этой теоремы, оператор $T^\natural(\gamma)$ ограниченно обратим
[\eqref{eq:3.99}]. Соответственно [\ref{20:sod}, \ref{prop:2.2}], оператор
$T^\bullet\rightleftharpoons (I^*)^{-1}TI^{-1}$, очевидным образом расширяющий
оператор $T^\circ$, является самосопряжённым.

Целью настоящего параграфа является разработка аналогичного метода построения
самосопряжённых расширений для симметрических операторных матриц описанного
в~\ref{pt:1.1} вида.

\subsection\label{pt:2.1}
Пусть в гильбертовом пространстве $\fH=\fH_1\oplus\fH_2$ задана симметрическая
операторная матрица \eqref{eq:1.99} с указанными в \ref{pt:1.1} свойствами.
Зафиксируем два значения $\varkappa,\,\tau\in \Bbb R$, удовлетворяющие соотношениям
$$
	\displaylines{(\forall y\in\dom T^\circ_{11})\qquad\langle(T_{11}^\circ-
		\varkappa)y,y\rangle_{\fH_1}\geqslant\|y\|_{\fH_1}^2,\cr
	(\forall y\in\dom T^\circ_{22})\qquad\langle(\tau-T^\circ_{22})y,
		y\rangle_{\fH_2}\geqslant\|y\|_{\fH_2}^2.}
$$
Обозначим через $\fD_2$ пополнение линейного множества $\dom T^\circ_{22}$
по норме
$$
	\|y\|_{\fD_2}\rightleftharpoons\langle(\tau-T_{22}^\circ)y,
		y\rangle_{\fH_2}^{1/2},
$$
через $T^\bullet_{22}$~--- связанное с пространством $\fD_2$ фридрихсовское
расширение оператора $T^\circ_{22}$, а через $\fD_1$~--- пополнение линейного
множества $\dom T_{11}^\circ$ по норме
$$
	\|y\|_{\fD_1}\rightleftharpoons\left[\langle(T_{11}^\circ-\varkappa)y,
		y\rangle_{\fH_1}+\langle(\tau-T^\bullet_{22})^{-1}T_{21}^\circ y,
		T_{21}^\circ y\rangle_{\fH_2}\right]^{1/2}.\leqno(\equation)
$$\label{eq:d1}%
По построению, при этом существуют обладающие плотными образами непрерывные
операторы вложения $I_1\colon\fD_1\to\fH_1$ и $I_2\colon\fD_2\to\fH_2$.
Это позволяет ввести в рассмотрение оснащение ${\bf A}\rightleftharpoons\{\fH,
\fD,\fD^*,I,I^*\}$, где положено $\fD\rightleftharpoons\fD_1\oplus\fD_2$
и $I\rightleftharpoons I_1\oplus I_2$, а через $\fD^*$ обозначено пространство,
сопряжённое пространству $\fD$.

Заметим, что при любом выборе вектора $y\in\fD_1$, удовлетворяющего дополнитель\-ному
условию $I_1y\in\dom T_{11}^\circ$, выполняются соотношения
$$
	\eqalignno{\|I_1^*(T_{11}^\circ-\varkappa)I_1 y\|_{\fD_1^*}&=
		\sup\limits_{\|z\|_{\fD_1}=1}|\langle(T_{11}^\circ-\varkappa)I_1 y,
		I_1 z\rangle_{\fH_1}|\cr
	&\leqslant\sup\limits_{\scriptstyle \|z\|_{\fD_1}=1\atop
		\scriptstyle I_1z\in\dom T_{11}^\circ}
		\langle(T_{11}^\circ-\varkappa)I_1 y,I_1 y\rangle_{\fH_1}^{1/2}
		\cdot\langle(T_{11}^\circ-\varkappa)I_1 z,
		I_1 z\rangle_{\fH_1}^{1/2}\cr
	&\leqslant\|y\|_{\fD_1},\cr
	\|I_2^*T_{21}^\circ I_1 y\|_{\fD_2^*}&=
		\sup\limits_{\|z\|_{\fD_2}=1}|\langle T_{21}^\circ I_1 y,
		I_2 z\rangle_{\fH_2}|\cr
	&=\sup\limits_{\|z\|_{\fD_2}=1}|\langle(\tau-T^\bullet_{22})^{-1/2}
		T_{21}^\circ I_1y,(\tau-T^\bullet_{22})^{1/2}
		I_2z\rangle_{\fH_2}|\cr
	&\leqslant\|(\tau-T^\bullet_{22})^{-1/2}T_{21}^\circ I_1y\|_{\fH_2}\cdot
		\sup\limits_{\|z\|_{\fD_2}=1} \|(\tau-T^\bullet_{22})^{1/2}
		I_2z\|_{\fH_2}\cr
	&\leqslant\|y\|_{\fD_1}.}
$$
Соответственно, замыканиями операторов $I_1^*T_{11}^\circ I_1$ и $I_2^*T_{21}^\circ
I_1$ являются некоторые всюду определённые ограниченные операторы $T_{11}\colon
\fD_1\to\fD_1^*$ и $T_{21}\colon\fD_1\to\fD_2^*$.

Аналогичным образом, при любом выборе вектора $y\in\fD_2$, удовлетворяющего
дополнительному условию $I_2y\in\dom T_{22}^\circ$, выполняются соотношения
$$
	\eqalignno{\|I_2^*(T_{22}^\circ-\tau)I_2y\|_{\fD_2^*}&=
		\sup\limits_{\|z\|_{\fD_2}=1}|\langle(\tau-T_{22}^\circ)I_2y,
		I_2 z\rangle_{\fH_2}|\cr
	&\leqslant\sup\limits_{\scriptstyle \|z\|_{\fD_2}=1\atop
		\scriptstyle I_2z\in\dom T_{22}^\circ}
		\langle(\tau-T_{22}^\circ)I_2 y,I_2 y\rangle_{\fH_2}^{1/2}
		\cdot\langle(\tau-T_{22}^\circ)I_2 z,I_2 z\rangle_{\fH_2}^{1/2}\cr
	&=\|y\|_{\fD_2},\cr
	\|I_1^*T_{12}^\circ I_2 y\|_{\fD_1^*}&=
		\sup\limits_{\scriptstyle \|z\|_{\fD_1}=1\atop
		\scriptstyle I_1 z\in\dom T_{11}^\circ}|\langle I_2 y,
		T_{21}^\circ I_1 z\rangle_{\fH_2}|\cr
	&=\sup\limits_{\scriptstyle \|z\|_{\fD_1}=1\atop
		\scriptstyle I_1 z\in\dom T_{11}^\circ}
		|\langle(\tau-T^\bullet_{22})^{1/2}I_2 y,
		(\tau-T^\bullet_{22})^{-1/2}T_{21}^\circ I_1 z\rangle_{\fH_2}|\cr
	&\leqslant\|(\tau-T^\bullet_{22})^{1/2}I_2y\|_{\fH_2}\cdot
		\sup\limits_{\scriptstyle \|z\|_{\fD_1}=1\atop
		\scriptstyle I_1 z\in\dom T_{11}^\circ}
		\|(\tau-T^\bullet_{22})^{-1/2}T_{21}^\circ I_1 z\|_{\fH_2}\cr
	&\leqslant\|y\|_{\fD_2}.}
$$
Соответственно, замыканиями операторов $I_2^*T_{22}^\circ I_2$ и $I_1^*T_{12}^\circ
I_2$ являются некоторые всюду определённые ограниченные операторы $T_{22}\colon
\fD_2\to\fD_2^*$ и $T_{12}\colon\fD_2\to\fD_1^*$.

Объединяя установленные факты, убеждаемся, что замыканием оператора $I^*T^\circ I$
является некоторый всюду определённый ограниченный оператор
$T\in {\cal B}({\bf A})$. Отвечающий ему оператор $T^\bullet\rightleftharpoons
(I^*)^{-1}TI^{-1}$ мы далее будем называть {\it угловым расширением} исходной
операторной матрицы $T^{\circ}$.

\subsubsection\label{prop:autoad}
{\it Угловое расширение $T^\bullet$ операторной матрицы $T^\circ$ представляет
собой самосопряжённый оператор в гильбертовом пространстве $\fH$.
}

\proof
Вещественнозначность квадратичной формы оператора $T^\bullet$ немедленно
вытекает из его определения. Соответственно, для доказательства предло\-жения
достаточно [\cite{Kato:1972}, Гл.~V, Теорема~4.3] установить
ограниченную обратимость оператора
$$
	T^\bullet-\pmatrix{\varkappa&0\cr 0&\tau}=(I^*)^{-1}\pmatrix{T_{11}-
		\varkappa I_1^*I_1&T_{12}\cr T_{21}&T_{22}-\tau I_2^*I_2}I^{-1}.
$$
Последняя, в свою очередь, немедленно вытекает из ограниченной обратимости
оператора
$$
	\pmatrix{T_{11}-\varkappa I_1^*I_1&T_{12}\cr T_{21}& T_{22}-
		\tau I_2^*I_2}=\pmatrix{1& T_{12}D_\tau^{-1}\cr 0&1}
		\pmatrix{S_{\varkappa,\tau}&0\cr 0&D_\tau}
		\pmatrix{1&0\cr D_\tau^{-1}T_{21}&1},\leqno(\equation)
$$\label{eq:FS}%
где положено $D_\tau\rightleftharpoons T_{22}-\tau I_2^*I_2$
и $S_{\varkappa,\tau}\rightleftharpoons T_{11}-\varkappa I_1^*I_1-T_{12}
D_\tau^{-1}T_{21}$: действительно, числовые области значений операторов $D_\tau$
и $S_{\varkappa,\tau}$ отделены от нуля, что означает [\ref{prop:St}]
ограниченную обратимость этих операторов.
\endproof

В ходе доказательства предложения \ref{prop:autoad} нами, по существу, была также
установлена справедливость следующего предложения.

\subsubsection\label{prop:anpr}
{\it Любое значение $\lambda\in\Bbb R$ принадлежит спектру оператора $T^\bullet$
в том и только том случае, когда оно принадлежит спектру линейного пучка
$T^\natural\colon\Bbb C\to{\cal B}({\bf A})$ вида
$$
	T^\natural(\lambda)=T-\lambda I^*I.
$$
При этом кратность любого собственного значения $\lambda\in\Bbb R$ оператора
$T^\bullet$ в точности совпадает с размерностью ядра оператора
$T^\natural(\lambda)$} [\ref{20:sod}, \ref{prop:2.2}, \ref{20:ker}].

\subsection
Связь разработанной выше в настоящем параграфе конструкции с результатами работы
[\cite{KLT}] даётся следующим легко проверяемым предложением.

\subsubsection\label{prop:KLT}
{\it Пусть матрица $T^\circ$ самосопряжена в существенном, операторы $T_{11}^\circ$
и $T_{22}^\circ$ полуограничены (соответственно, снизу и сверху), а также выполнено
хотя бы одно из следующих условий:

\smallskip
\itemitem{$1^\circ$.} Оператор $T_{22}^\circ$ самосопряжён в существенном, причём
для некоторых $\gamma,\,\varkappa\in\Bbb R$ выполняется соотношение
$$
	(\forall y\in\dom T_{11}^\circ)\qquad \|T_{21}^\circ y\|_{\fH_1}^2
		\leqslant\gamma\cdot\langle(T_{11}^\circ-\varkappa)y,
		y\rangle_{\fH_1}.
$$

\itemitem{$2^\circ$.} Для некоторых $\gamma,\varkappa,\tau\in\Bbb R$ выполняются
соотношения
$$
	\displaylines{(\forall y\in\dom T_{11}^\circ)\qquad \|T_{21}^\circ
		y\|_{\fH_1}^2\leqslant\gamma\cdot\langle(T_{11}^\circ-
		\varkappa)y,y\rangle_{\fH_1},\cr
	(\forall y\in\dom T_{22}^\circ)\qquad \|T_{12}^\circ y\|_{\fH_2}^2
		\leqslant\gamma\cdot\langle(\tau-T_{22}^\circ)y,
		y\rangle_{\fH_2}.}
$$

\itemitem{$3^\circ$.} Операторы $T_{11}^\circ$ и $T_{22}^\circ$ являются ограниченными.

\smallskip\noindent
Тогда замыкание оператора $T^\circ$ является его угловым расширением.
}

\medskip
Фигурирующие в формулировке предложения \ref{prop:KLT} три случая взаимных
соотношений между элементами операторной матрицы $T^\circ$ суть в точности
три типа взаимного доминирования этих элементов, рассмотренные в работе
[\cite{KLT}].

\subsection
В качестве простого примера применения процедуры углового расширения рассмотрим
действующую в пространстве $L_2[0,1]\oplus L_2[0,1]$ симметрическую операторную
матрицу
$$
	\pmatrix{-d^2/dx^2&-d^2/dx^2\cr -d^2/dx^2&0}\leqno(\equation)
$$\label{eq:10ekz}%
с областью определения $\Wo_2^2[0,1]\oplus\Wo_2^2[0,1]$. Этой матрице отвечают
пространства $\fD_1=\Wo_2^2[0,1]$ и $\fD_2=L_2[0,1]$ (с некоторыми эквивалентными
стандартным нормами), а само расширение задаётся той же матрицей с областью
определения $\Wo_2^2[0,1]\oplus W_2^2[0,1]$. При этом операторы $T_{22}-
\lambda I_2^*I_2$ представляют собой действующие в $L_2[0,1]$ операторы умножения
на постоянную $\lambda$, что означает в случае $\lambda\neq 0$ равносильность
ограниченной обратимости оператора $T^\natural(\lambda)=T-\lambda I^*I$
и передаточного оператора
$$
	S(\lambda)\rightleftharpoons T_{11}-\lambda I_1^*I_1+\lambda^{-1}T_{12}
		T_{21}\colon\Wo_2^2[0,1]\to\Wo_2^{-2}[0,1].
$$
Последний, в свою очередь, имеет вид
$$
	\langle S(\lambda)y,z\rangle\equiv\int_0^1\bigl[y'\overline{z'}-
		\lambda y\overline{z}+\lambda^{-1}y''\overline{z''}\bigr]\,dx,
$$
что означает равносильность его ограниченной обратимости отсутствию нетривиальных
решений граничной задачи
$$
	\displaylines{y^{(4)}-\lambda y''-\lambda^2 y=0,\cr
		y(0)=y'(0)=y(1)=y'(1)=0.}
$$
Согласно предложению \ref{prop:anpr}, спектр последней граничной задачи в точности
совпадает на множестве $\Bbb R\setminus\{0\}$ со спектром определённой
на $\Wo_2^2[0,1]\oplus W_2^2[0,1]$ операторной матрицы \eqref{eq:10ekz}.

Отметим, что матрица \eqref{eq:10ekz} не относится ни к одному из трёх типов,
рассмотренных в работе [\cite{KLT}] и указанных в формулировке предложения
\ref{prop:KLT}.


\section{Вариационные принципы}\label{par:3}
\subsection
Предложение \ref{prop:anpr} сводит вопрос о спектре углового расширения $T^\bullet$
операторной матрицы $T^\circ$ к вопросу о спектре линейного пучка $T^\natural$
ограниченных операторов вида
$$
	T^\natural(\lambda)=\pmatrix{T_{11}-\lambda I_1^*I_1&T_{12}\cr
		T_{21}&T_{22}-\lambda I_2^*I_2}.\leqno(\equation)
$$\label{eq:Tnat}%
В случае, когда значение $\lambda\in\Bbb C$ принадлежит резольвентному множеству
оператора $T_{22}^\bullet$, оператор $T^\natural_{22}(\lambda)=T_{22}-\lambda I_2^*
I_2$ обладает ограниченным обратным [\ref{prop:2.2}], что означает возможность
факторизации операторной матрицы \eqref{eq:Tnat} в виде \eqref{eq:FS}.
Соответственно, вне спектра оператора $T^\bullet_{22}$ собственные значения
оператора $T^\bullet$ совпадают~--- с сохранением кратностей~--- с собственными
значениями оператор-функции $S$ вида
$$
	S(\lambda)\rightleftharpoons T_{11}-\lambda I_1^*I_1-T_{12}
		[T^\natural_{22}(\lambda)]^{-1}T_{21}.
$$
Настоящий параграф будет посвящён описанию основанных на таком замечании
вариационных принципов для собственных значений оператора $T^\bullet$.

\subsection\label{pt:3.1}
Далее через $\ind A$ всегда обозначается отрицательный индекс инерции
квадратич\-ной формы оператора $A$, т.~е. точная верхняя грань размерностей
подпространств $\frak M\subseteq\dom A$, удовлетворяющих условию
$$
	(\exists\varepsilon>0)\,(\forall y\in\frak M\setminus\{0\})\qquad
		\langle Ay,y\rangle<0.
$$
Имеет место следующий факт.

\subsubsection\label{prop:vp1}
{\it Пусть отрезок $[\zeta^-,\zeta^+]\subseteq\varrho(T_{22}^\bullet)$ таков,
что оператор $S(\zeta^+)\colon\fD_1\to\fD_1^*$ представляет собой вполне непрерывное
возмущение некоторого равномерно положи\-тельного оператора. Тогда спектр оператора
$T^\bullet$ на полуинтервале $[\zeta^-,\zeta^+)$ чисто дискретен, причём его
суммарная кратность равна величине $\ind S(\zeta^+)-\ind S(\zeta^-)$.
}

\proof
Зафиксируем число $\varepsilon>0$, представляющее собой нижнюю оценку некоторого
вполне непрерывного возмущения оператора $S(\zeta^+)$. Ввиду справедли\-вости для любых
чисел $\lambda_1\in[\zeta^-,\zeta^+)$, $\lambda_2\in(\lambda_1,\zeta^+]$ и вектора
$y\in\fD_1$ соотношений
$$
	\leqalignno{\langle [S(\lambda_2)-S(\lambda_1)]y,y\rangle&=
		\int_{\lambda_1}^{\lambda_2}{d\langle S(\lambda)y,y\rangle\over
		d\lambda}\,d\lambda&\cr
	&=\int_{\lambda_1}^{\lambda_2}\left(-\|I_1y\|_{\fH_1}^2-
		\|I_2[T^\natural_{22}(\lambda)]^{-1}T_{21}y\|_{\fH_2}^2\right)\,
		d\lambda\cr
	&\leqslant-(\lambda_2-\lambda_1)\,\|I_1y\|^2_{\fH_1},&
		(\equation)}
$$\label{eq:sootn3}%
постоянная $\varepsilon$ будет минорировать некоторое вполне непрерывное возмущение
оператора $S(\lambda)$ при любом выборе значения $\lambda\in [\zeta^-,\zeta^+]$.
Положим теперь
$$
	\Lambda_n(\lambda)\rightleftharpoons\inf_{\scriptstyle\frak M\subseteq\fD_1
		\atop\scriptstyle\dim\frak M>n}\sup_{\scriptstyle
		y\in\frak M\atop\scriptstyle\|y\|_{\fD_1}=1}
		\inf\Bigl\{\langle S(\lambda)y,y\rangle,\varepsilon\Bigr\}.
$$
Непрерывный и монотонный характер зависимости операторов $S(\lambda)$ от параметра
гарантирует непрерывность и невозрастание каждой из функций $\Lambda_n\colon
[\zeta^-,\zeta^+]\to\Bbb R$, где $n\in\Bbb N$.

Стандартные вариационные принципы для ограниченных самосопряжённых опера\-торов
[\cite{AG:1966}, п.~82, Теорема~2] утверждают, что независимо от выбора значения
$\lambda\in[\zeta^-,\zeta^+]$ размерность подпространства $\ker S(\lambda)$ равна
величине
$$
	\#\{n\in\Bbb N\;:\;\Lambda_n(\lambda)=0\},
$$
причём в случае обращения этой размерности в нуль оператор $S(\lambda)$ обладает
ограниченным обратным. Кроме того, из оценок \eqref{eq:sootn3} вытекает факт
отрицательной определённости квадратичной формы оператора $S(\lambda_2)$ на любом
подпространстве, на котором неположительна квадратичная форма оператора
$S(\lambda_1)$. Тем самым, каждая из функций $\Lambda_n$ строго убывает в своих
нулях. Соответственно, пересечение спектра оператора $T^\bullet$ с полуинтервалом
$[\zeta^-,\zeta^+)$ представляет собой [\ref{prop:anpr}] дискретное множество
всевозможных нулей функций $\Lambda_n$, причём суммарная кратность этой части
спектра оператора $T^\bullet$ совпадает с величиной
$$
	\#\{n\in\Bbb N\;:\; (\Lambda_n(\zeta^-)\geqslant 0)\mathbin{\hbox{\&}}
		(\Lambda_n(\zeta^+)<0)\}=\ind S(\zeta^+)-\ind S(\zeta^-).
	\leqno(\equation)
$$\label{eq:val}%
Тем самым, доказываемое предложение является верным.
\endproof

Проведённые при доказательстве предложения \ref{prop:vp1} рассуждения в основном
воспроизводят доказательство теоремы [\cite{Vl:2003}, Теорема~2]. Аналогичным
образом доказыва\-ется также следующее предложение.

\subsubsection\label{prop:vp2}
{\it Пусть отрезок $[\zeta^-,\zeta^+]\subseteq\varrho(T_{22}^\bullet)$ таков,
что оператор $S(\zeta^-)\colon\fD_1\to\fD_1^*$ представляет собой вполне непрерывное
возмущение некоторого равномерно отрица\-тельного оператора. Тогда спектр оператора
$T^\bullet$ на полуинтервале $(\zeta^-,\zeta^+]$ чисто дискретен, причём его
суммарная кратность равна величине $\ind [-S(\zeta^-)]-\ind [-S(\zeta^+)]$.
}

\subsection
Имеют место следующие четыре факта.

\subsubsection\label{prop:3.10}
{\it Пусть $\zeta\in\varrho(T_{22}^\bullet)$~--- имеющее конечную кратность
изолированное собственное значение оператора $T^\bullet$. Пусть также
для некоторого $\zeta^-<\zeta$ справедливо соотношение $[\zeta^-,\zeta)\subseteq
\varrho(T^\bullet)\cap\varrho(T_{22}^\bullet)$, причём оператор $S(\zeta^-)$
представляет собой вполне непрерывное возмущение некоторого равномерно
положительного оператора. Тогда найдётся такое $\zeta^+>\zeta$, что справедливо
соотношение $(\zeta,\zeta^+]\subseteq\varrho(T^\bullet)\cap
\varrho(T_{22}^\bullet)$, а оператор $S(\zeta^+)$ представляет собой вполне
непрерывное возмущение некоторого равномерно положительного оператора.
}

\proof
В рассматриваемом случае найдётся [\cite{Tren:1980}, п.~21.4] самосопряжён\-ный
оператор конечного ранга $K\colon\fH\to\fH$, для которого оператор $T^\bullet+
K-\zeta$ будет ограниченно обратимым. Без умаления общности рассмотрения можно
считать оператор $K$ настолько малым по норме, что оператор
$$
	\tilde D_\zeta\rightleftharpoons T_{22}+I_2^*(K_{22}-\zeta)I_2
$$
обладает ограниченным обратным. В этом случае определённая на некотором содержащем
точку $\zeta$ интервале $\Gamma\subseteq(\zeta^-,+\infty)$ оператор-функция
$\tilde S$ вида
$$
	\tilde S(\lambda)\rightleftharpoons T_{11}+I_1^*(K_{11}-\lambda)I_1-
		(T_{12}+I_1^*K_{12}I_2)\,\tilde D_\lambda^{-1}\cdot
		(T_{21}+I_2^*K_{21}I_1)
$$
принимает при $\lambda=\zeta$ ограниченно обратимое значение [\ref{prop:anpr}].

Ввиду непрерывности функции $\tilde S$ в смысле равномерной операторной топологии,
найдётся содержащий точку $\zeta$ интервал $\Delta\subseteq\Gamma$ со следующими
свойствами:

\smallskip
\itemitem{$1^\circ$.} При любом выборе значения $\lambda\in\Delta$ оператор
$\tilde S(\lambda)$ обладает ограниченным обратным.

\itemitem{$2^\circ$.} Операторы $\theta(\tilde S(\lambda))$, где через $\theta$
обозначена функция Хэвисайда, непрерывно зависят от параметра $\lambda\in\Delta$
в смысле равномерной операторной топологии.

\smallskip
Зафиксируем теперь произвольное расположенное слева от точки $\zeta$ значение
$\mu\in\Delta$. Оператор $\tilde S(\mu)$ представляет собой имеющее конечный ранг
возмущение оператора $S(\mu)$, по условию являющегося вполне непрерывным
возмущением некоторого равномерно положительного оператора. Соответственно,
оператор $\tilde S(\mu)$ сам является вполне непрерывным возмущением некоторого
равномерно положительного оператора, и потому оператор $\theta(\tilde S(\mu))$
представляет собой проектор на подпространство конечной коразмерности. Согласно
вышесказанному, проекторами на некоторые подпространства той же конечной
коразмерности [\cite{AG:1966}, п.~39] будут выступать также операторы
$\theta(\tilde S(\lambda))$ при всех $\lambda\in\Delta$. Последнее как раз
и означает, что при любом выборе значения $\lambda\in\Delta$ операторы
$\tilde S(\lambda)$ и $S(\lambda)$ представляют собой вполне непрерывные
возмущения некоторого равномерно положительного оператора.
\endproof

\subsubsection\label{prop:3.11}
{\it Пусть $\zeta\in\varrho(T_{22}^\bullet)$~--- имеющее конечную кратность
изолированное собственное значение оператора $T^\bullet$. Пусть также
для некоторого $\zeta^+>\zeta$ справедливо соотношение $(\zeta,\zeta^+]\subseteq
\varrho(T^\bullet)\cap\varrho(T_{22}^\bullet)$, причём оператор $S(\zeta^+)$
представляет собой вполне непрерывное возмущение некоторого равномерно
отрицательного оператора. Тогда найдётся такое $\zeta^-<\zeta$, что справедливо
соотношение $[\zeta^-,\zeta)\subseteq\varrho(T^\bullet)\cap
\varrho(T_{22}^\bullet)$, а оператор $S(\zeta^-)$ представляет собой вполне
непрерывное возмущение некоторого равномерно отрицательного оператора.
}

\medskip
Доказательство предложения \ref{prop:3.11} полностью аналогично доказательству
предложе\-ния \ref{prop:3.10}.

\subsubsection
{\it Пусть все расположенные на отрезке $[\zeta^-,\zeta^+]\subseteq
\varrho(T_{22}^\bullet)$ точки спектра оператора $T^\bullet$ являются
изолированными собственными значениями конечной крат\-ности. Пусть также оператор
$S(\zeta^-)$ представляет собой вполне непрерывное возмущение некоторого равномерно
положительного оператора. Тогда оператор $S(\zeta^+)$ также представляет собой
вполне непрерывное возмущение некоторого равномерно положительного оператора.
}

\subsubsection
{\it Пусть все расположенные на отрезке $[\zeta^-,\zeta^+]\subseteq
\varrho(T_{22}^\bullet)$ точки спектра оператора $T^\bullet$ являются
изолированными собственными значениями конечной крат\-ности. Пусть также оператор
$S(\zeta^+)$ представляет собой вполне непрерывное возмущение некоторого равномерно
отрицательного оператора. Тогда оператор $S(\zeta^-)$ также представляет собой
вполне непрерывное возмущение некоторого равномерно отрицательного оператора.
}

\medskip
Два последних предложения легко выводятся из предложений \ref{prop:3.10}
и \ref{prop:3.11} индукцией по числу лежащих на отрезке $[\zeta^-,\zeta^+]$
собственных значений оператора $T^\bullet$.

\subsection
Связь результатов настоящего параграфа с результатами работы [\cite{KLT}]
устанав\-ливается следующим образом.

Обозначим через $\lambda_e$ точную нижнюю грань пересечения существенного
спектра оператора $T^\bullet$ с некоторой полупрямой $(\zeta,+\infty)\subseteq
\varrho(T_{22}^\bullet)$, где $\zeta\in\varrho(T^\bullet)$. В качестве точной
нижней грани пустого множества здесь и далее понимается величина $+\infty$.
Обозначим также через $\{\lambda_n\}_{n=0}^{\infty}$ последовательность,
членами которой выступают либо $(n+1)$-ые снизу (с учётом кратности) собственные
значения оператора $T^\bullet$ из интервала $(\zeta,\lambda_e)$, либо величина
$\lambda_e$, если требуемое собственное значение не существует. Имеет место
следующий тривиальным образом вытекающий из результатов настоящего параграфа факт.

\subsubsection\label{prop:VP}
{\it Пусть величина $\ind S(\zeta)$ конечна. Тогда при любом $n\in\Bbb N$
выполняется равенство
$$
	\lambda_n=\inf\{\lambda\in (\zeta,+\infty)\;:\;\ind S(\lambda)>
		\ind S(\zeta)+n\}.
$$
}

Основные результаты работы [\cite{KLT}] получаются комбинированием предложений
\ref{prop:VP} и \ref{prop:KLT}.


\section{Дополнительные иллюстрации}\label{par:4}
\subsection\label{pt:4.1}
Несложно заметить, что для получения основных результатов пункта~\ref{pt:2.1}
полуограниченность оператора $T_{11}^\circ$ в действительности не является
необходимой. Достаточным оказывается выполнение следующих двух условий:

\smallskip
\itemitem{$1^\circ$.} Существуют величины $\varkappa,\tau\in\Bbb R$, для которых
квадратичная форма
$$
	\frak s[y]\rightleftharpoons\langle(T_{11}^\circ-\varkappa)
		y,y\rangle_{\fH_1}+\langle(\tau-T_{22}^\bullet)^{-1}
		T_{21}^\circ y,T_{21}^\circ y\rangle_{\fH_2}
$$
равномерно по $\fH_1$-норме знакоопределена на линейном множестве
$\dom T_{11}^\circ$.

\itemitem{$2^\circ$.} При выборе в качестве пространства $\fD_1$ пополнения линейного
множества $\dom T_{11}^\circ$ по норме $\bigl|\frak s[y]\bigr|^{1/2}$ замыканием
оператора $I^*T^{\circ}I$ является некоторый всюду определённый ограниченный оператор
$T\in{\cal B}({\bf A})$.

\smallskip
В качестве примера здесь может быть рассмотрена действующая в пространстве
$L_2[0,1]\oplus L_2[0,1]$ самосопряжённая операторная матрица
$$
	\pmatrix{id/dx&-d/dx\cr d/dx&0}
$$
с областью определения ${\tilde W}_2^1[0,1]\oplus{\tilde W}_2^1[0,1]$,
где снова использовано обозначение \eqref{eq:Wtild}. На роль пространств $\fD_1$
и $\fD_2$ могут быть взяты пространства ${\tilde W}_2^1[0,1]$ и $L_2[0,1]$,
соответственно. Аналогичными проведённым в пункте \ref{pt:2.1} рассуждениями
легко устанавливается, что сосредоточенный на $\varrho(T_{22}^\bullet)=
\Bbb C\setminus\{0\}$ спектр рассматриваемой матрицы совпадает со спектром
граничной задачи
$$
	\displaylines{-y''+i\lambda y'-\lambda^2=0,\cr y(0)-y(1)=y'(0)-y'(1)=0.}
$$

\subsection
Результаты предыдущего параграфа были сформулированы применительно к задаче
о спектре простейшего линейного пучка $T^\bullet-\lambda$. Между тем, как несложно
заметить, они допускают естественное распространение на случай самосопряжённой
оператор-функции достаточно общего вида. А именно, пусть $T\colon[\zeta^-,
\zeta^+]\to{\cal B}\,({\bf A})$ есть гладкая в смысле сильной операторной топологии
оператор-функция, для которой при любом выборе значения $\lambda\in [\zeta^-,
\zeta^+]$ соответствующий оператор $[T(\lambda)]^\bullet$ самосопряжён.
Пусть при этом также выполнены следующие два условия:

\smallskip
\itemitem{$1^\circ$.} Существует вещественное число $\varepsilon>0$, для которого
независимо от выбора значения $\lambda\in [\zeta^-,\zeta^+]$ правая нижняя
компонента $T_{22}(\lambda)$ операторной матрицы $T(\lambda)$ удовлетворяет
соотношению
$$
	(\forall y\in\fD_2)\qquad\langle T_{22}(\lambda)y,y\rangle
		\leqslant-\varepsilon\,\|y\|_{\fD_2}^2.
$$

\itemitem{$2^\circ$.} Определённая на отрезке $[\zeta^-,\zeta^+]$ оператор-функция
$S$ вида
$$
	S(\lambda)\rightleftharpoons T_{11}(\lambda)-T_{12}(\lambda)
		[T_{22}(\lambda)]^{-1}T_{21}(\lambda)
$$
непрерывна в смысле равномерной операторной топологии, а её значения суть
вполне непрерывные возмущения некоторых равномерно положительных операторов.

\smallskip
Во избежание недоразумений подчеркнём, что норма пространств $\fD_1$ и $\fD_2$
здесь уже не предполагается каким-либо образом связанной с соотношениями
из пункта~\ref{pt:2.1}.

Напомним [\cite{LSY}, \cite{Vl:2003}], что собственное значение $\mu\in (\zeta^-,
\zeta^+)$ самосопряжённой оператор-функции $\lambda\mapsto [T(\lambda)]^\bullet$
считается {\it имеющим отрицательный тип}, если для любого ненулевого вектора
$y\in\ker[T(\mu)]^\bullet$ числовая функция $\lambda\mapsto
\langle[T(\lambda)]^\bullet y,y\rangle_{\fH}$ определена на некотором содержащем
точку $\mu$ интервале и имеет в этой точке строго отрицательную производную.
При сделанных предположениях имеют место следующие два факта.

\subsubsection
{\it Каждое имеющее отрицательный тип изолированное конечнократное собст\-венное
значение $\mu\in(\zeta^-,\zeta^+)$ оператор-функции $\lambda\mapsto
[T(\lambda)]^\bullet$ является имеющим отрицательный тип изолированным
собственным значением той же кратности для передаточной оператор-функции $S$.
}

\proof
При сделанных предположениях из теоремы Банаха--Штейн\-га\-уза
[\cite{Tren:1980}, п.~11.5] вытекает равномерная по параметру $\lambda\in
[\zeta^-,\zeta^+]$ ограниченность операторов $T(\lambda)$ и связанная с этим
возможность дифференцирования оператор-функции $S$ в смысле сильной операторной
топологии. Согласно предложению \ref{prop:anpr} и замечаниям из пункта~\ref{pt:4.1},
точка $\mu$ является изолирован\-ным собственным значением оператор-функции $S$,
причём соответствующее собственное подпространство допускает [\ref{20:ker}]
представление
$$
	\left\{y\in\fD_1\;:\;I\pmatrix{y\cr z_{\mu}}\in
		\ker [T(\mu)]^\bullet\right\},
$$
где использовано сокращение $z_{\lambda}\rightleftharpoons -[T_{22}(\lambda)]^{-1}
T_{21}(\lambda)y$. Учёт выполняющихся для любого вектора $y\in\ker S(\mu)$ равенств
$$
	\postdisplaypenalty 20000
	\eqalignno{\left.{d\langle S(\lambda)y,y\rangle\over d\lambda}
		\right|_{\lambda=\mu}&=\left.{d\langle S(\lambda)y,
		y\rangle\over d\lambda}\right|_{\lambda=\mu}-
		2\Re\left\langle\pmatrix{S(\mu)& 0\cr 0&T_{22}(\mu)}
		\pmatrix{ 0\cr (dz_{\lambda}/d\lambda)|_{\lambda=\mu}},
		\pmatrix{y\cr 0}\right\rangle\cr
	&=\left.{d\over d\lambda}\left\langle\pmatrix{S(\lambda)&
		0\cr 0&T_{22}(\lambda)}\pmatrix{y\cr z_{\mu}-z_{\lambda}},
		\pmatrix{y\cr z_{\mu}-z_{\lambda}}\right\rangle
		\right|_{\lambda=\mu}\cr
	&=\left.{d\over d\lambda}\left\langle [T(\lambda)]^\bullet I
		\pmatrix{y\cr z_{\mu}}, I\pmatrix{y\cr z_{\mu}}\right\rangle_{\fH}
		\right|_{\lambda=\mu}}
$$
завершает доказательство.
\endproof

\subsubsection\label{prop:3.4.1}
{\it Пусть все расположенные на интервале $(\zeta^-,\zeta^+)$ собственные значения
оператор-функции $T$ имеют отрицательный тип, а каждый из операторов $T(\zeta^\pm)$
обладает ограниченным обратным. Тогда спектр оператор-функции $\lambda\mapsto
[T(\lambda)]^\bullet$ на интервале $(\zeta^-,\zeta^+)$ чисто дискретен, причём
его суммарная кратность равна величине $\ind S(\zeta^+)-\ind S(\zeta^-)$.
}

\proof
Зафиксируем произвольную точку $\mu\in\Bbb R$ и рассмотрим определённую
на некотором содержащем эту точку интервале гладкую оператор-функцию $F$ вида
$$
	F(\lambda)\rightleftharpoons\pmatrix{F_{11}(\lambda)&
		F_{12}(\lambda)\cr F_{21}(\lambda)&F_{22}(\lambda)},
$$
значениями которой выступают ограниченные самосопряжённые операторы в ортого\-нальной
прямой сумме $\frak E_1\oplus\frak E_2$ некоторых конечномерных
гильбертовых пространств. Предположим также выполнение равенств $F_{11}(\mu)=0$,
$F_{21}(\mu)=0$ и обратимость операторов $F_{11}'(\mu)$ и $F_{22}(\mu)$.
Тогда при $\lambda\to\mu+0$ оказываются справедливыми асимптотики
$$
	\eqalign{\ind F(\lambda)&=\ind F_{11}(\lambda)+\ind [F_{22}(\lambda)-
		F_{21}(\lambda)F_{11}^{-1}(\lambda)F_{12}(\lambda)]\cr
	&=\ind [(\lambda-\mu)\cdot (F_{11}'(\mu)+o(1))]+\ind [F_{22}(\mu)
		+O(\lambda-\mu)]\cr
	&=\ind F_{11}'(\mu)+\ind F(\mu).}
$$

Повторим теперь в основном рассуждения из доказательства предложения
\ref{prop:vp1}. А именно, сопоставим оператор-функции $S$ последовательность
непрерывных числовых функций $\Lambda_n\colon[\zeta^-,\zeta^+]\to\Bbb R$. Применяя
полученную асимптотику индексов инерции к случаю, когда пространство
$\frak E_1$ представляет собой ядро оператора $S(\mu)$, пространство
$\frak E_2$~--- отвечающее отрицательной части спектра инвариантное
подпространство того же оператора, а квадратичные формы операторов $F(\lambda)$
совпадают с таковыми для операторов $S(\lambda)$, убеждаемся, что каждая
из функций $\Lambda_n$ строго убывает в своих нулях. Согласно предложению
\ref{prop:anpr} и замечаниям из пункта~\ref{pt:4.1}, это означает совпадение
суммарной кратности расположенного на интервале $(\zeta^-,\zeta^+)$ спектра
оператор-функции $\lambda\mapsto [T(\lambda)]^\bullet$ с величиной \eqref{eq:val}.
\endproof


\vskip 0.4cm
\eightpoint\rm
{\leftskip 0cm\rightskip 0cm plus 1fill\parindent 0cm
\bf Литература\par\penalty 20000}\vskip 0.4cm\penalty 20000
\bibitem{KLT} M.~Kraus, M.~Langer, C.~Tretter. {\it Variational principles
and eigenvalue estimates for unbounded block operator matrices and applications}~//
Journ. of Comput. and Appl. Math.~--- 2004.~--- V.~171.~--- P.~311--334.
\bibitem{LM} Ж.-Л.~Лионс, Э.~Мадженес. {\it Неоднородные граничные задачи
и их приложения}. М.: Мир, 1971.
\bibitem{BSh:1983} Ф.$\,$А.~Березин, М.$\,$А.~Шубин. {\it Уравнение Шрёдингера}.
М.: Изд-во МГУ, 1983.
\bibitem{Kato:1972} Т.~Като. {\it Теория возмущений линейных операторов}.
М.: Мир, 1972.
\bibitem{AG:1966} Н.$\,$И.~Ахиезер, И.$\,$М.~Глазман. {\it Теория линейных
операторов в гильбертовом пространстве}. М.: Наука, 1966.
\bibitem{Vl:2003} А.$\,$А.~Владимиров. {\it Оценки числа собственных значений
самосопряжённых оператор-функций}~// Матем. заметки.~--- 2003.~--- Т.~74,
\No~6.~--- С.~838--847.
\bibitem{Tren:1980} В.$\,$А.~Треногин. {\it Функциональный анализ}. М.: Наука, 1980.
\bibitem{LSY} P.~Lancaster, A.~Shkalikov, Qiang~Ye. {\it Strongly definitizable
linear pencils in Hilbert space}~// Integr. Equat. Oper. Th.~--- 1993.~---
V.~17.~--- P.~338--360.
\bibitem{LLT} H.~Langer, M.~Langer, C.~Tretter. {\it Variational principles
for eigenvalues of block operator matrices}~// Indiana Univ. Math. Jour.~---
2002.~--- V.~51, \No~6.~--- P.~1427--1459.
\end